\pdfoutput=1
\RequirePackage{ifpdf}
\ifpdf 
\documentclass[pdftex]{sigma}
\else
\documentclass{sigma}
\fi

\numberwithin{equation}{section}

\newtheorem{Theorem}{Theorem}[section]
\newtheorem*{Theorem*}{Theorem}

\newtheorem{Question}[Theorem]{Question}

\theoremstyle{definition}

\newcommand{\E}{{\mathbb{E}}}
\newcommand{\Vol}{{\rm Vol}}
\newcommand{\B}{{\mathbb{B}}}
\newcommand{\conv}{{\rm conv}}

\begin{document}
\allowdisplaybreaks

\newcommand{\arXivNumber}{2406.18428}

\renewcommand{\PaperNumber}{109}

\FirstPageHeading

\ShortArticleName{Small Volume Bodies of Constant Width with Tetrahedral Symmetries}

\ArticleName{Small Volume Bodies of Constant Width\\ with Tetrahedral Symmetries}

\Author{Andrii ARMAN~$^{\rm a}$, Andriy BONDARENKO~$^{\rm b}$, Andriy PRYMAK~$^{\rm a}$ and Danylo RADCHENKO~$^{\rm c}$}

\AuthorNameForHeading{A.~Arman, A.~Bondarenko, A.~Prymak and D.~Radchenko}

\Address{$^{\rm a)}$~Department of Mathematics, University of Manitoba, Winnipeg, MB, R3T 2N2, Canada}
\EmailD{\mail{andrew0arman@gmail.com}, \mail{prymak@gmail.com}}
\URLaddressD{\url{https://sites.google.com/view/a-arman/home}, \url{http://prymak.net/}}

\Address{$^{\rm b)}$~Department of Mathematical Sciences, Norwegian University of Science and Technology,\\
\hphantom{$^{\rm b)}$}~NO-7491 Trondheim, Norway}
\EmailD{\mail{andriybond@gmail.com}}
\URLaddressD{\url{https://www.ntnu.edu/employees/andrii.bondarenko}}

\Address{$^{\rm c)}$~Universit\'{e} de Lille, CNRS, Laboratoire Paul Painlev\'{e}, F-59655 Villeneuve d'Ascq, France}
\EmailD{\mail{danradchenko@gmail.com}}
\URLaddressD{\url{https://danrad.net/}}

\ArticleDates{Received June 04, 2025, in final form December 06, 2025; Published online December 21, 2025}	

\Abstract{For every $n\ge 2$, we construct a body $U_n$ of constant width $2$ in $\mathbb{E}^n$ with small volume and symmetries of a regular $n$-simplex. $U_2$ is the Reuleaux triangle. To the best of our knowledge, $U_3$ was not previously constructed, and its volume is smaller than the volume of other three-dimensional bodies of constant width with tetrahedral symmetries. While the volume of $U_3$ is slightly larger than the volume of Meissner's bodies of width $2$, it exceeds the latter by less than $0.137\%$. For all large $n$, the volume of $U_n$ is smaller than the volume of the ball of radius $0.891$.}


\Classification{52A20; 52A15; 52A23; 52A40; 28A75; 49Q20}

\section{Introduction}\label{sec:intro}

A convex body in the $n$-dimensional Euclidean space $\E^n$ is a convex compact set with non-empty interior. A convex body has constant width $w$ if its orthogonal projection onto any line has length~$w$. For a comprehensive overview of bodies of constant width and related topics, see~\cite{Maetal}.\looseness=-1

The largest by volume convex body of constant width $2$ in $\E^n$ is the unit ball $\B^n$ in $\E^n$, which is a consequence of the Urysohn's inequality (see, e.g., \cite[equation~(7.21), p.~382]{Schn}). Blaschke~\cite{Bl} and Lebesgue~\cite{Le} showed that smallest area bodies of constant width in $\E^2$ are the Reuleaux triangles. It is an open question (Blaschke--Lebesgue problem) to find the smallest possible volume of a convex body of a fixed constant width in $\E^n$, $n\ge 3$, and to characterize the minimizers. For $n=3$, Meissner and Schilling~\cite{MeSc} showed how to modify the direct generalization of the Reuleaux triangles to the three dimensional case to obtain two non-congruent bodies of constant width which are now usually referred to as Meissner's bodies. It was conjectured by Bonnesen and Fenchel~\cite{BoFe} that Meissner's bodies are the minimum volume bodies of constant width in $\E^3$. A reader interested in Meissner's bodies is referred to the work of Kawohl and Weber~\cite{Ka-We}. We would like to credit the survey~\cite{Ho} of Horv\'ath for the historical details and references.

For the asymptotic case when the dimension $n\to\infty$, it has been an open question raised by Schramm~\cite{Schr} (see also the survey of Kalai~\cite[Problem~3.4]{K}) if there exist convex bodies of constant width having the volume significantly smaller than the volume of the ball of the same width. We have recently answered this question in the affirmative in a joint work with Nazarov~\cite{ABNPR}. Namely, we constructed the sets
\[
M_n:=\bigl\{v-w\mid v,w\in\E^n_+,\, |v|^2+\bigl(|w|+\sqrt{2}\bigr)^2\le4 \bigr\},
\]
where $\E^n_+$ is the subset of $\E^n$ consisting of the points whose all coordinates are non-negative and~$|\cdot|$ is the Euclidean norm. It is proved in~\cite{ABNPR} that $M_n$ are convex bodies of width $2$ and that~${\Vol(M_n)<\Vol(0.891 \B^n)}$ for large $n$.

We say that a convex body $U \subset \E^n$ has the symmetries of a regular $n$-simplex (the convex polyhedron in $\E^n$ with $n+1$ equidistant vertices) if there exists a regular $n$-simplex $\Delta$ such that the full symmetry group of $\Delta$ acts as symmetries of $U$; equivalently,
$
\tau(U) = U $ for every orthogonal transformation $\tau$ satisfying $\tau(\Delta) = \Delta$.

Define $U_n$ as the orthogonal projection of $M_{n+1}$ onto the $n$-dimensional subspace of $\E^{n+1}$ orthogonal to the diagonal direction $(1,\dots,1)\in\E^{n+1}$, i.e.,
\begin{equation*}
	U_{n}=\operatorname{proj}_{(1,\ldots,1)^{\perp}} (M_{n+1})=\{\operatorname{proj}_{(1,\ldots,1)^{\perp}} (v) \mid v \in M_{n+1}\}.
\end{equation*}

Then $U_n$ is also a body of constant width $2$. Moreover, as $M_{n+1}$ is invariant under any permutation of coordinates, $U_n$ has the symmetries of a regular $n$-simplex. We obtained the following result.

\begin{Theorem}
	$U_n$ is a body of constant width $2$ in $\E^n$ with the symmetries of a regular $n$-simplex.
\end{Theorem}

Now let us discuss the volume of $U_n$. It is easily checked that $U_2$ is a Reuleaux triangle, so it has the smallest area among all planar bodies of constant width. For $U_3$, we have the following statement.

\begin{Theorem}\label{prop:vol}
	We have $\frac{\Vol(U_3)}{\Vol(\B^3)} = 0.8029702551499\dots $, and moreover
	\begin{align*}
	\Vol(U_3) = {}& \frac{4\pi}{3}\bigl(4-\sqrt{6}\bigr)+24\left(2-\sqrt{2}-\arctan\frac1{\sqrt{8}}\right)\\
&-96\int_{0}^{\frac1{\sqrt{3}}}\sqrt{\frac{1+a^2}{1-a^2}}\arctan\left(\frac{a}{\sqrt{1+a^{2}}+1}\right){\rm d}a .
	\end{align*}
\end{Theorem}
The proof can be found in~Appendix~\ref{sec:appA}. For comparison, the volume of Meissner's bodies of width $2$ is (see, e.g., \cite[p.~68]{CG})
\[
8\pi \left(\frac23-\frac{\sqrt{3}}{4}\cdot\arccos\frac13 \right)= \Vol\bigl(\B^3\bigr) \cdot 0.80187362\dots\, .
\]
We see that $U_3$ is larger in volume than Meissner's bodies, but only by less than $0.137\%$, while~$U_3$ has symmetries of a regular tetrahedron ($3$-simplex). The smallest previously known body of constant width having such symmetries appears to be the Minkowski's average of the two Meissner's bodies. Its volume exceeds the volume of the Meissner's bodies by slightly more than~$0.2\%$, as stated in~\cite[p.~4]{Ka-We2} referring to a work of Oudet which appears to be no longer available online, but is consistent with our computations as outlined in Appendix~\ref{sec:appB}.

Another example of a three-dimensional body of constant width with tetrahedral symmetries is Robert's body~\cite{Ro}, which is a representative of the class of peabodies~\cite{pea}. The volume of the Robert's body exceed the volume of the Meissner's bodies by more than $0.23\%$.

\begin{figure}[t]
	\centering
	\includegraphics[width=0.33\linewidth]{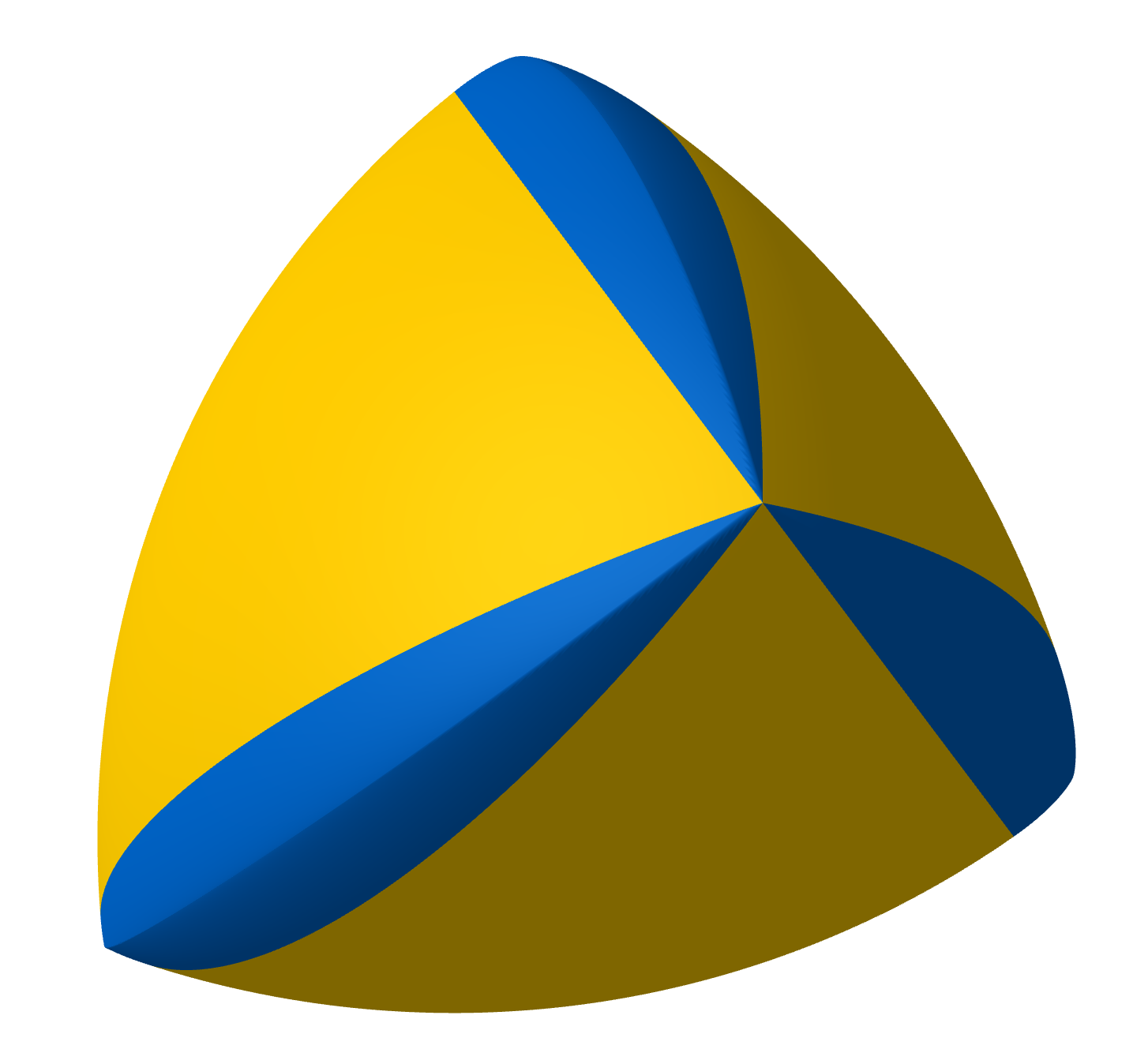}
	\qquad
	\includegraphics[width=0.33\linewidth]{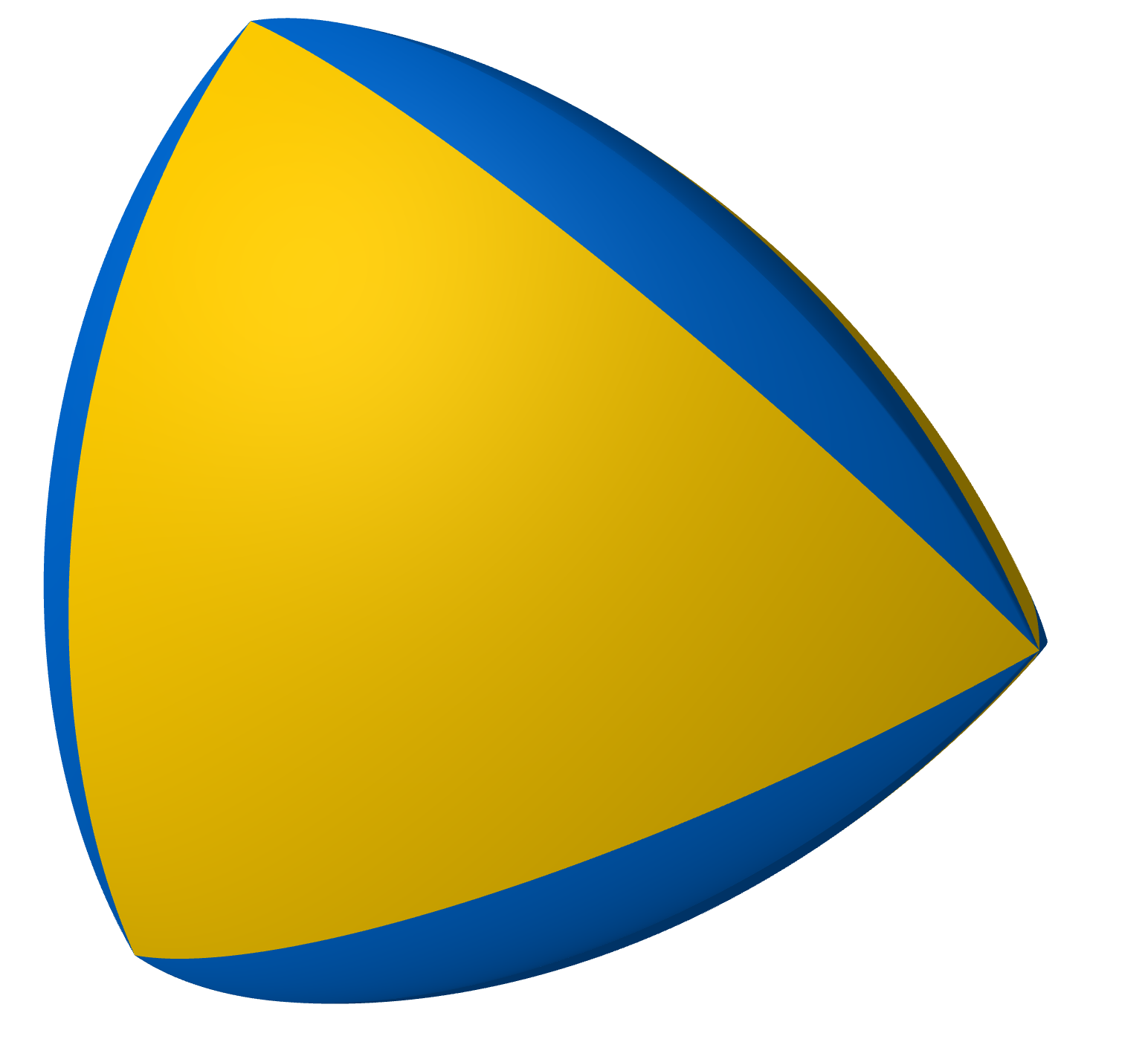}
	\caption{Two views of $U_3$.}
	\label{fig:u3}
\end{figure}

In~Figure~\ref{fig:u3}, we provide two views of $U_3$.\footnote{A three-dimensional model is available at \url{https://www.danrad.net/misc/ubody.html}.} The smooth pieces of the boundary of $U_3$ consist of the yellow pieces of spheres from the Reuleaux tetrahedron, and two blue surfaces per each edge, which properly cut the Reuleaux tetrahedron to make it a body of constant width. Our construction of $U_3$ is based on an orthogonal projection of a $4$-dimensional body $M_4$ defined explicitly through Cartesian coordinates. In contrast, Meissner's bodies are constructed by a~specific three-dimensional procedure of rounding certain edges of the Reuleaux tetrahedron using rotation of an appropriate circular arc about an edge. The rounding is performed on three edges adjacent to a vertex to get one of the Meissner's bodies, and on three opposite edges to get the other Meissner's body, then resulting Minkowski's average of the two will have the symmetries of $3$-simplex. An interested reader may find more details of constructions of Meissner's bodies in, e.g.,~\cite[Section~8.3, p.~171]{Maetal} or~\cite{Ka-We2}. We also remark that the four spherical parts/``faces'' of~$U_3$ (yellow in~Figure~\ref{fig:u3}) and the Meissner's bodies and their average coincide, the difference is near the edges.

In contrast to the conjecture of Bonnesen and Fenchel, according to~\cite[p.~68]{CG}, \cite[p.~34, Problem~A22]{CFG} and~\cite[p.~261]{GSch} it was conjectured by Danzer that the minimizer of the volume among bodies of constant width must have the symmetry group of a regular simplex. The example of $U_3$ puts Meissner's bodies as better candidates for the minimizers in $\E^3$.

A necessary condition for a three-dimensional body of constant width to be a local minimizer of the volume was obtained by Anciaux and Guilfoyle~\cite{AG}. Namely, the smooth components of the minimizer have their smaller principle curvature constant and equal to the reciprocal of the width of the body. The Meissner's bodies as well as $U_3$ do satisfy this condition.

\begin{Question}
	Is $U_n$ the minimizer of the volume among all bodies of constant width in $\E^n$ having symmetries of a regular simplex?
\end{Question}
Recall that $U_2$ is a Reuleaux triangle, so the answer for $n=2$ is affirmative. Even for $n=3$, it would be very interesting to learn the answer. For small $n>3$, it would be interesting to compare $\Vol(U_n)$ to the volume of the bodies obtained by applying the raising dimension process of Lachand-Robert and Oudet~\cite{LaOu} to the Meissner's bodies.

We conclude by considering the higher dimensional case, where the problem can be reduced to the estimate of the volume of $M_{n+1}$.
\begin{Theorem}
	There exists positive integer $n_0$ such that
	$\Vol(U_n)<\Vol(0.891\B^n)$ for all $n\ge n_0$.
\end{Theorem}
\begin{proof}
	We claim that for any unit $u\in\E^{n+1}$, the projection $W':=\operatorname{proj}_{u^\perp}(W)$ of a convex body~$W$ of constant width $2$ in $\E^{n+1}$ onto the subspace $u^\perp=\{x\mid x\cdot u=0\}$ satisfies
	\begin{equation}\label{eq:1}
		\Vol_{n}(W')\le \frac{n+1}2 \Vol_{n+1}(W).
	\end{equation}
Indeed, if $S:=S_{u^\perp}W$ is the Steiner symmetral of $W$ with respect to $u^\perp$ (see, e.g.,~\cite[Section~10.3, p.~536]{Schn}), then the section of $S$ by $u^\perp$ is precisely $W'$ and the width of $S$ in the direction of $u$ is exactly $2$ with some two points $v_{\pm}\in S$ satisfying $v_{\pm}\cdot u=\pm1$ (here we use the fact that in a~body of constant width the points of contact with any two parallel supporting hyperplanes form a segment perpendicular to these hyperplanes, and are exactly the width apart, see, e.g.,~\cite[Section~3.1, Theorem~3.1.1, p.~59]{Maetal}). By convexity, each $C_{\pm}:=\conv(W'\cup \{v_{\pm}\})$ is a subset of $S$, so by $\Vol_{n+1}(C_{\pm})=\frac1{n+1} \Vol_{n}(W')$ and $\Vol_{n+1}(C_+)+\Vol_{n+1}(C_-)\le \Vol_{n+1}(S)=\Vol_{n+1}(W)$, hence (\ref{eq:1}) follows.
	
	Now the proposition follows from $\Vol_{n+1}(M_{n+1})<\Vol_{n+1}\bigl(0.8908\B^{n+1}\bigr)$ (see~\cite[Remark~1]{ABNPR}) and the fact that \smash{$\bigl(\frac{\Vol_{n+1}(\B^{n+1})}{\Vol_n(\B^n)}\bigr)^{1/n}\to1$}, $n\to\infty$.
\end{proof}

We remark that $0.891$ can be replaced with any constant bigger than the value described in~\cite[Remark~1]{ABNPR}. Also note that $\Vol_n(U_n)$ cannot be much smaller than $\Vol_{n+1}(M_{n+1})$ since the width of $M_{n+1}$ in the direction of $(1,1,\dots,1)$ is $2$, so $2\Vol_n(U_n)\ge \Vol_{n+1}(M_{n+1})$.

\appendix

\section[Volume computation for U\_3]{Volume computation for $\boldsymbol{U_3}$}\label{sec:appA}
\begin{figure}[h]
	\centering
	\includegraphics[width=0.4\linewidth]{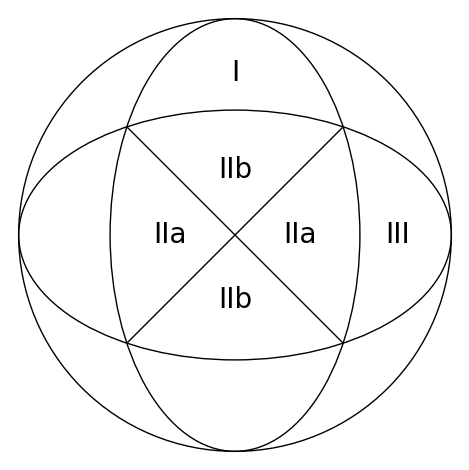}
	\caption{Integration regions in $(a,b)$ coordinates. The ellipses are given by equations $a^2+3b^2=1$ and $3a^2+b^2=1$.}
	\label{fig:sch}
\end{figure}

By~\cite[Theorem~2]{AG}, the volume of a three-dimensional body of constant width 2 is given in terms of its support function by
\begin{equation}\label{eq:vol}
	\frac{4\pi}{3}-\int_{\mathbb{S}^2}\frac12|\nabla h(\theta)|^2-(h(\theta)-1)^2{\rm d}\mu(\theta) .
\end{equation}
The support function of $M_n$, for all $n$, is given by
\begin{equation*} \label{eq:supp}
	h(\theta) = \begin{cases}
		\sqrt{2}|\theta_{+}| ,&|\theta_{+}|\ge |\theta_{-}| ,\\
		2-\sqrt{2}|\theta_{-}| ,& |\theta_{+}|\le |\theta_{-}| ,
	\end{cases}
\end{equation*}
where $\theta$ is decomposed in a unique way as $\theta=\theta_+-\theta_-$, with $\theta_+,\theta_-\in\E^n_+$ and $\theta_+\cdot\theta_-=0$.
The support function of $U_3$ is the restriction of $h$ to the 3-plane $\{w+x+y+z=0\}$.
Up to symmetries, we need to integrate~\eqref{eq:vol} in 3 different types of regions: $w,x,y>0>z$ (I), $w,x>0>y,z$ (II), and $w>0>x,y,z$ (III). The case II splits into two more cases according to whether $w^2+x^2\ge y^2+z^2$ (IIa) or $w^2+x^2< y^2+z^2$ (IIb). The support function of $M_4$ can then be computed as
\[
h(w,x,y,z) = \begin{cases}
	2+\sqrt{2}z ,&\mbox{case I} ,\\
	\sqrt{2\bigl(w^2+x^2\bigr)} ,& \mbox{case II}_a ,\\
	2-\sqrt{2\bigl(y^2+z^2\bigr)} ,& \mbox{case II}_b ,\\
	\sqrt{2}w ,& \mbox{case III} .
\end{cases}
\]
We pick an orthonormal basis in $\{w+x+y+z=0\}$ consisting of the vectors
\begin{equation*}
	e_1= \left(\frac{1}{\sqrt{2}},-\frac{1}{\sqrt{2}},0,0\right) ,\qquad
	e_2= \left(0,0,\frac{1}{\sqrt{2}},-\frac{1}{\sqrt{2}}\right) ,\qquad
	e_3= \left(\frac{1}{2},\frac{1}{2},-\frac{1}{2},-\frac{1}{2}\right) .
\end{equation*}
Denoting the coordinates in the basis $e_1$, $e_2$, $e_3$ by $(a,b,c)$, we have the following relations:
\begin{gather*}
a=\frac{w-x}{\sqrt{2}} ,\qquad
b=\frac{y-z}{\sqrt{2}} ,\qquad
c=w+x=-y-z .
\\
w=\frac{c+a\sqrt{2}}{2} ,\qquad
x=\frac{c-a\sqrt{2}}{2} ,\qquad
y=\frac{-c+b\sqrt{2}}{2} ,\qquad
z=\frac{-c-b\sqrt{2}}{2} .
\end{gather*}
Case I then corresponds to $b\sqrt{2}\ge c\ge |a|\sqrt{2}$, Case~IIa corresponds to $c\ge |a|\sqrt{2}\ge |b|\sqrt{2}$, Case~IIb corresponds to $c\ge |b|\sqrt{2}\ge |a|\sqrt{2}$, and Case III corresponds to $a\sqrt{2}\ge c\ge |b|\sqrt{2}$. Since~${a^2+b^2+c^2=1}$, we can rewrite these regions in terms of just $a$ and $b$. Figure~\ref{fig:sch} shows the different regions in $(a,b)$ coordinates.

We can then write the support function of $U_3$ in new coordinates as follows:
\[
h_{U_3}(a,b)-1 = \begin{cases}
	1-\frac{c}{\sqrt{2}}-b ,& \mbox{case I} ,\\
	\sqrt{1+a^2-b^2}-1 ,& \mbox{case II}_a ,\\
	1-\sqrt{1-a^2+b^2} ,& \mbox{case II}_b ,\\
	\frac{c}{\sqrt{2}}+a-1 ,&\mbox{case III} ,
\end{cases}
\]
where $c=\sqrt{1-a^2-b^2}$.
Note that the area form on the sphere in $(a,b)$ coordinates is \smash{$\frac{{\rm d}a{\rm d}b}{\sqrt{1-a^2-b^2}}$}. Because of the symmetry $a\leftrightarrow b$ it suffices to deal with the cases IIa and III. The spherical gradient $|\nabla h|^2$ can be computed from Euclidean gradient by projecting to the tangent space of~$\mathbb{S}^2$, so in Case III we have
\[\frac12|\nabla h|^2-(h-1)^2=\frac34-\frac12\bigl(a+c/\sqrt{2}\bigr)^2-\bigl(a+c/\sqrt{2}-1\bigr)^2 ,\]
and in Case IIa we have
\[\frac12|\nabla h|^2-(h-1)^2=\frac12\left(\frac{a^2+b^2}{1+a^2-b^2}-\frac{(a^2-b^2)^2}{1+a^2-b^2}\right)-\bigl(\sqrt{1+a^2-b^2}-1\bigr)^2 .\]
Collecting the contributions to the integral in~\eqref{eq:vol} according to symmetries, we get
\[\Vol(U_3) = \frac{4\pi}{3} -16I_1-48I_2 ,\]
where
\[I_1 = \int_{0}^{1/2}\int_{\sqrt{\frac{1-b^2}{3}}}^{\sqrt{1-3b^2}}
\left(\frac34-\frac12\bigl(a+c/\sqrt{2}\bigr)^2-\bigl(a+c/\sqrt{2}-1\bigr)^2\right)\frac{{\rm d}a}{\sqrt{1-a^2-b^2}}{\rm d}b\]
and
\[I_2 = \int_{0}^{1/2}\int_{b}^{\sqrt{\frac{1-b^2}{3}}}
\left(\frac{a^2}{1+a^2-b^2}-\frac{3\bigl(a^2-b^2\bigr)}{2}-2+2\sqrt{1+a^2-b^2}\right)\frac{{\rm d}a}{\sqrt{1-a^2-b^2}}{\rm d}b .\]
The integral defining $I_1$ can be computed symbolically to obtain
\[I_1 = \frac{\pi\sqrt{6}-3\sqrt{2}}{12}-\arctan\bigl(\sqrt{2}/5\bigr) .\]
After some manipulations, $I_2$ can be rearranged into
\[I_2 = \frac{7\sqrt{2}}{12}-1-\arctan\bigl(\sqrt{2}\bigr)+\pi/4+2\int_{0}^{1/2}\int_{b}^{\sqrt{\frac{1-b^2}{3}}}
\frac{\sqrt{1+a^2-b^2}}{\sqrt{1-a^2-b^2}}{\rm d}a{\rm d}b .\]
Combining these and simplifying, we get
\[\Vol(U_3) = \frac{4\pi}{3}\bigl(4-\sqrt{6}\bigr)+24\bigl(2-\sqrt{2}\bigr)-96\int_{0}^{1/2}\int_{b}^{\sqrt{\frac{1-b^2}{3}}}\frac{\sqrt{1+a^2-b^2}}{\sqrt{1-a^2-b^2}}{\rm d}a{\rm d}b .\]
The last integral can be further rewritten as
\begin{gather*}
	\int_{0}^{1/2}\int_{b}^{\sqrt{\frac{1-b^2}{3}}}\frac{\sqrt{1+a^2-b^2}}{\sqrt{1-a^2-b^2}}{\rm d}a{\rm d}b\\
\qquad
	= \int_{0}^{1/\sqrt{3}}\sqrt{\frac{1+a^2}{1-a^2}}\left(\frac{a}{2\bigl(1+a^2\bigr)}+\arctan\left(\frac{a}{\sqrt{1+a^{2}}+1}\right)\right){\rm d}a \\
\qquad	=\frac{\arctan\bigl(1/\sqrt{8}\bigr)}{4}+\int_{0}^{1/\sqrt{3}}\sqrt{\frac{1+a^2}{1-a^2}}\arctan\left(\frac{a}{\sqrt{1+a^{2}}+1}\right){\rm d}a ,
\end{gather*}
but it does not seem to reduce to any more elementary expression. The above integral can then be computed to any desired precision using, say, PARI/GP, giving the numerical value
\[\frac{\Vol(U_3)}{\Vol\bigl(\B^3\bigr)} = 0.80297025514991011046814277\dots\, .\]

\section[Volume computation of Minkowski's average of Meissner's bodies]{Volume computation of Minkowski's average\\ of Meissner's bodies}\label{sec:appB}

We follow the approach from Appendix~\ref{sec:appA}. Let $A$ be a Meissner's body with vertices \smash{$\bigl(\pm1,0,\frac{1}{\sqrt{2}}\bigr)$} and \smash{$\bigl(0,\pm1,-\frac1{\sqrt{2}}\bigr)$}, in $(a,b,c)$ coordinates, whose edge joining \smash{$\bigl(\pm1,\frac{1}{\sqrt{2}},0\bigr)$} is not rounded, i.e., coinciding with the Reuleaux tetrahedron with the same vertices. Let $B$ be the second Meissner's body with the same vertices and opposite edges rounded. The support function of either $A$ or~$B$ is the same as the support function of $U_3$ for case I and case III, which leads to the same term~$-16I_1$ in $\Vol(U_3)$ as in $\Vol(A)$, $\Vol(B)$ and $\Vol\bigl(\frac{A+B}{2}\bigr)$. Now it suffices to calculate the support function near the non-rounded edge of the Reuleaux tetrahedron which will give the support function of the opposite rounded edge of the Meissner's body by the constant width property. We have
\[
h_A(a,b)-1=\begin{cases}
	1-\frac{c}{\sqrt{2}}-b, & 0\le a\le \frac12,\quad \frac12\le b\le \sqrt{\frac{1-a^2}3},\\[1mm]
	\sqrt{3\bigl(1-b^2\bigr)}-\frac{c}{\sqrt{2}}-1, &0\le b\le \frac12,\quad 0\le a\le \sqrt{\frac{1-b^2}3}.
\end{cases}
\]
By symmetry and the definitions of $A$ and $B$, $h_B(a,b)+h_A(b,a)=2$, so we arrive at
\[
h_{\frac{A+B}{2}}(a,b)-1=\begin{cases}
	1-\frac b2-\frac12 \sqrt{3\bigl(1-a^2\bigr)}, & 0\le a\le \frac12,\quad \frac12\le b\le \sqrt{\frac{1-a^2}3},\\[1mm]
	\frac12 \sqrt{3\bigl(1-b^2\bigr)}-\frac12 \sqrt{3\bigl(1-a^2\bigr)}, & 0\le a\le b\le \frac12.
\end{cases}
\]
Now using~\eqref{eq:vol}, computing the spherical gradients and applying available symmetries, we obtain
\[
\Vol\left(\frac{A+B}{2}\right)=\frac{4\pi}{3}-16I_1-48(I_3+I_4),
\]
where
\begin{align*}
I_3={}&\frac18\int_0^{\frac12}\int_{\frac12}^{\sqrt{\frac{1-a^2}3}}
\left(
\frac{1+2a^2}{1-a^2}-2\big(b-2+\sqrt{3\big(1-a^2\big)}\big)^2-\left(\frac{\sqrt{3}a^2}{\sqrt{1-a^2}}-b\right)^2
\right)\\
&\times\frac{{\rm d}b}{\sqrt{1-a^2-b^2}} {\rm d}a,
\end{align*}
and
\begin{align*}
	I_4={}&\frac18\int_0^{\frac12}\int_a^{\frac12}
	\left(\frac{3a^2}{1-a^2}+\frac{3b^2}{1-b^2} - 6\bigl(\sqrt{1-b^2}-\sqrt{1-a^2}\bigr)^2 \right.\\
& \left.- 3\left(\frac{a^2}{\sqrt{1-a^2}}-\frac{b^2}{\sqrt{1-b^2}}\right)^2
	\right)\frac{{\rm d}b}{\sqrt{1-a^2-b^2}} {\rm d}a.
\end{align*}
Numerical integration gives
\[\frac{\Vol\bigl(\frac{A+B}{2}\bigr)}{\Vol\bigl(\B^3\bigr)} = 0.803806345386\dots \,.\]

\subsection*{Acknowledgements}

We would like to thank anonymous referees for carefully reading the paper and providing valuable feedback. A.~Arman acknowledges support in part by a postdoctoral fellowship of the Pacific Institute for the Mathematical Sciences. A.~Bondarenko was supported in part by Grant 334466 of the Research Council of Norway, and A.~Prymak was supported by NSERC of Canada Discovery Grant RGPIN-2020-05357. D.~Radchenko acknowledges funding by the European Union (ERC, FourIntExP, 101078782).

\pdfbookmark[1]{References}{ref}
\LastPageEnding

\end{document}